\documentclass[12pt,leqno]{amsart}

\usepackage{amssymb}
\usepackage{latexsym}

\newcommand\dspace{\lineskip=2pt\baselineskip=18pt\lineskiplimit=0pt}
\dspace

\newcounter{deficislo}[section]
\newcounter{lemacislo}[section]
\newcounter{propcislo}[section]
\newcommand{\lema}[1]{\refstepcounter{lemacislo}{\noindent \bf
Lemma \thesection.\thelemacislo.\  }{\it #1}}
\newcommand{\defi}[1]{\refstepcounter{deficislo}{\noindent \bf
Definition \thesection.\thedeficislo.\  }{\it #1}}

\newcommand{\prop}[1]{\refstepcounter{propcislo}{\noindent \bf
Theorem \thesection.\thepropcislo.\  }{\it #1}}

\newcommand{\pz}{P(z,\bar z)}

\newcommand{\ep}{\epsilon}
\newcommand{\be}{\beta}
\newcommand{\ga}{\gamma}
\newcommand{\de}{\delta}

\newcommand{\zz}{(z,\bar z)}

\newcommand{ \al}{\alpha}

\newcommand{\la}{\lambda}
\newcommand{\La}{\Lambda}

\begin{document}
\title{
The Catlin  multitype
and biholomorphic equivalence of models}
\author {Martin Kol\'a\v r}
\address {Department of Mathematics and Statistics, Masaryk University,
Kotlarska 2,  611 ~37 Brno } 
\email {mkolar@math.muni.cz }
\maketitle

\begin{abstract}
We consider an alternative 
approach to a fundamental
CR invariant -- the Catlin multitype.
It is applied to a general smooth hypersurface in $\mathbb
C^{n+1}$,  not necessarily pseudoconvex. Using this approach, we
prove biholomorphic equivalence of models, and give an explicit
description of biholomorphisms between different models. A constructive finite
algorithm for computing the multitype is described. 
The
results can be
viewed as providing a necessary step in understanding local
biholomorphic equivalence of Levi degenerate hypersurfaces
of finite Catlin multitype.

\end{abstract}

\section{Introduction}
The subject of this paper is  local  biholomorphic geometry of
Levi-degenerate  hypersurfaces in $\mathbb C^{n+1}$,  and a fundamental CR
invariant -- the Catlin multitype. We consider a constructive
approach, which allows to understand the local
equivalence problem on the level of weighted homogeneous models.

The problem of local biholomorphic equivalence for real hypersurfaces in  complex
space has a long history (we refer to the survey articles \cite{BER2}, \cite{We} for a historical account).
In recent years, the problem has been intensively studied on Levi degenerate manifolds, 
mostly using the extrinsic approach of Poincar\' e and 
Moser. 
In fact, a result of 
\cite{KZ} indicates that the intrinsic approach of  Cartan, Chern and Tanaka is in general 
not available in the degenerate setting.

We start by 
 reviewing  some motivating facts from complex dimension two.
The lowest order CR invariant of a
smooth hypersurface $M \subseteq \mathbb C^2$ at a point $p \in M$ is the
type of the point, introduced by J. J. Kohn in  his pioneering work \cite{K}.
The type is an integer measuring the maximal order of
contact between $M$ and complex curves passing through $p$.
 In terms of coordinates,
 the point
 is of finite type $k$ if and only if there exist local holomorphic 
coordinates $(z,w)$ such that
the defining equation for $M$ takes form
\begin{equation}
Im \;w = P(z, \bar z) + o(Re\; w,
\vert z \vert^k), \label{pzz}
\end{equation}
 where
$P$ is a nonzero  homogeneous polynomial of degree
$k$ without harmonic terms.
The manifold $ Im\; w = P(z,\bar z) $ is the
model hypersurface at $p$.
Here $P$ is determined uniquely up to  linear
transformations in the  complex tangential variable $z$, and 
 one  immediately obtains important invariants
from the coefficients of $P$ (see e.g. \cite{KN}, \cite{Ko1}).

To study higher order invariants,  consider 
a biholomorphic transformation
\begin{equation}
w^* = g(z,w),  \ \ \ \ \ \ \ \ \  z^*  = f(z,w),\label{ffgg}
\end{equation} 
which preserves the local description (\ref{pzz}). 
The main tool for analyzing the action of (\ref{ffgg}) on the defining equation of $M$
is  the  generalized Chern-Moser
operator
\begin{equation}
L(f,g) =  \mbox{Re}  \left\lbrace ig(z,Re\; w +i\pz) +
2 \frac{\partial P}{\partial z} f(z,Re\; w + i\pz)\right\rbrace \  ,
\end{equation}
whose existence
is a fundamental consequences of the finite type condition.
Examining
the kernel and image  of $L$ one can construct a  complete set of local invariants 
 (\cite{Ko1}).

In higher dimensions,
local geometry of Levi degenerate
hypersurfaces is far  more complicated, even on the initial level.
 Invariants relevant  for
analysis of the inhomogeneous Cauchy-Riemann equations
are now  obtained by considering orders of contact with singular complex varieties.
If $d_k$ denotes the maximal order of contact of $M$ with complex
varieties of dimension $k$ at $p$, the n-tuple $(d_n,\dots, d_1)$
is called the D'Angelo multitype of  $M$ at $p$ (\cite{D}).

For pseudoconvex hypersurfaces, D. Catlin (\cite{C}) introduced
a different  notion of multitype,  using a more algebraic
approach. The entries of the Catlin multitype take rational values, but 
need not be integers,
anymore.
This approach provides a
defining equation analogous to (\ref{pzz}), and a well defined
weighted-homogeneous model,
an essential
tool 
for local analysis
(see e.g. \cite{KK}, \cite{KY}).

There is a class of
 hypersurfaces on which the two multitypes coincide
(termed semiregular  \cite{DH}, 
or h-extendible 
\cite{Y}),
but in the most interesting instances, the two
multitypes are not equal.

 In this paper we use Catlin's definition of
multitype for a general smooth hypersurface in $\mathbb C^{n+1}$.
The definition itself is nonconstructive, and the
corresponding models are not uniquely defined.
In order to study higher order CR invariants
it becomes essential to understand the non-uniqueness in the 
definition  of models.
In particular, it is not a priori clear if all models have to be
biholomorphically equivalent (for pseudoconvex  h-extendible hypersurfaces this problem 
was considered in \cite{Nik}).

Hypersurfaces of finite Catlin multitype provide the natural class of manifolds
for which a generalization of the Chern-Moser operator is well defined.

 We denote 
again by $P$ the leading weighted homogeneous polynomial determined by Catlin's
construction, and 
consider a biholomorphic transformation

\begin{equation}
 w^* = w + g(z, w),\ \ \ \ \ z_i^* = z_i
+ f_i(z, w).
\end{equation}
The operator now takes form 
\begin{equation}
L(f,g) =   Re \left\lbrace ig(z, Re\; w+iP(z,\bar z)) +
2\sum_{j=1}^{n} \frac{\partial P}{\partial z_j}  f_j(z,Re\; w + i\pz)\right\rbrace.
\end{equation}
The first necessary step in understanding this operator is to consider
the strictly subhomogeneous level, in the sense of Definition 2.3 below.  
Our results imply, in particular, that the kernel of $L$ is always trivial on this level.
Analysis of the kernel and image of $L$, and applications to the local equivalence
problem is the subject of a forthcoming article.

The paper is organized as follows. In Section 2 we define the
Catlin multitype of a general smooth hypersurface in $\mathbb
C^{n+1}$. This leads to distinguished weighted coordinate systems.
Then  we consider the associated weighted homogeneous
transformations, and define their  subhomogeneous and
superhomogeneous analogs. In Section 3 we analyze model
hypersurfaces, and define a normalization, which is used
in an essential way in the following section.

Section 4 considers the  biholomorphic equivalence problem for
models. We prove that all models at a given point  are
biholomorphically equivalent, by explicitly described polynomial
transformations.
 Using this result
we give in Section 5 a constructive finite algorithm for computing the
multitype.

\section{Hypersurfaces of finite multitype}
Let  $M \subseteq \mathbb C^{n+1}$ be a smooth hypersurface (not
necessarily pseudoconvex),  and $p $ be a Levi degenerate point on
$M$. We will assume that $p$ is a point of finite type in the
sense of Bloom and Graham.
Throughout the paper, the standard multiindex notation will be
used.

 Let
$(z,w)$ be local holomorphic coordinates centered at $p$,
where $z =(z_1, z_2, ..., z_n)$ and  $z_j = x_j + iy_j$,
$w=u+iv$. The hyperplane $\{ v=0 \}$ is assumed to be tangent to
$M$ at $p$. We describe $M$  near $p$ as the graph of a uniquely
determined real valued function
\begin{equation} v = \Psi(z_1,\dots, z_n,  \bar z_1,\dots,\bar z_n,  u).
\label{vp}
\end{equation}

The definition of multitype is based on  weighted coordinate systems. Roughly
speaking, the weights measure the order of vanishing of a suitably
chosen  defining function in each of the variables.
As the first step, the weights of the complex nontangential
variables $w$, $u$ and $v$ are set equal to one. Then we consider
the complex tangential variables.
\\[2mm]
\defi{A weight is an n-tuple of nonnegative
 rational numbers $\La = (\la_1, ...,
\la_n)$, where $0 \leq\la_j\leq \frac12$, and $\la_j \ge
\la_{j+1}$, such that for each $k$ either $\la_k =0$, or  there exist nonnegative
integers $a_1, \dots, a_k$ satisfying $a_k > 0$ and
$$ \sum_{j=1}^k a_j \la_j = 1.$$
}
\\[2mm]

If $\La$ is a weight,   the weighted degree of a monomial $c_{\al
\beta l}z^{\al}\bar z^\beta u^{l} $ is defined to be
$$ 
l +  \sum_{i=1}^n (\al_i + \beta_i ) \la_i.$$
A polynomial $P(z, \bar z, u)$  is $\La$-homogeneous
of weighted degree $\kappa$ if it is a sum of
 monomials of weighted degree $\kappa$.
\\[2mm]
The weighted length of a multiindex $\al = (\al_1, \dots, \al_n)$ is
defined by

$$\vert \alpha\vert_{\La} = \la_1 \al_1 + \dots + \la_n \al_n.$$

Similarly, if $\al = (\al_1, \dots, \al_n)$ and  $\hat \al =
(\hat \al_1, \dots, \hat \al_n)$ are two multiindices, the weighted
length of the  pair $(\al, \hat \al)$ is
$$\vert (\alpha,\hat \al) \vert_{\La} = \la_1 (\al_1 +\hat \al_1) \dots +
\la_n (\al_n + \hat \al_n).$$
The weighted order of a differential operator
 $\dfrac{\partial^{\vert \al + \hat \al \vert + l}}{\partial z^{\al}\partial \bar z^{\hat \al} \partial u^l}$ 
is equal to $l + \vert (\alpha,\hat \al) \vert_{\La}$.

 A weight $\La$ will be called distinguished if there exist
local holomorphic coordinates $(z,w)$ in which the defining equation of $M$ takes form
\begin{equation} v = P\zz + o_{\La}(1),\label{1}\end{equation}
where $P\zz$ is a nonzero $\La$-homogeneous polynomial of
weighted degree one without pluriharmonic terms, and $o_{\La}(1)$
denotes a smooth function whose derivatives of weighted order less than or equal to 
one vanish.  

The fact that distinguished weights do exist follows from the
assumption of Bloom-Graham finite type (\cite{BG}).
\\[2mm]

\defi{Let  $\Lambda_M = (\mu_1, \dots, \mu_n)$  be the infimum of
distinguished weights with respect to the lexicographic ordering.
The multitype of $M$ at $p$ is defined to be the n-tuple $(m_1,
m_2, \dots, m_n)$, where $m_j = \frac1{\mu_j}$ if $\mu_j \neq 0$
and $m_j = \infty $ if $\mu_j = 0$. If none of the $m_j$ is
infinity, we say that $M$ is of finite multitype at $p$.
}
\\[2mm]
Note that since the definition of multitype  considers all
distinguished weights, the infimum is a biholomorphic invariant,
and we may speak of {\it the} multitype.

Coordinates corresponding to a distinguished weight $\Lambda$, in
which the local description of $M$ has form (\ref{1}), with $P$
being  $\Lambda$-homogeneous,  will be called $\La$-adapted.

$\Lambda_M $ will be called the multitype weight, and $\La_M$-adapted 
coordinates will be also referred to as multitype
coordinates.

 It is easy to verify that for
any $\delta > 0$ there are only finitely many possible rational 
values for any weight entry, satisfying $\la_i > \delta$. 
It follows  that if $M$ is of
finite multitype at $p$, $\La_M$ - adapted coordinates do exist.

From now on we assume that $p \in M$ is a point of finite
multitype.

Let $t$ denote the number of different entries appearing in the
multitype weight, and $\nu_j, j = 1, \dots, t,$ be the length of
the $j$-th constant piece of the multitype weight. Hence, denoting
$k_j = \sum_{i=1}^j \nu_i$,
we have
$$\mu_1 = ... =\mu_{k_1}> \mu_{k_1 + 1} = \dots = \mu_{k_2}
 > ...
 =  \mu_{k_{t-1}}> \mu_{k_{t-1}+1}=\dots = \mu_{n}. $$

We define a  'generating' sequence of weights $\La_1, ..., \La_t$
as follows. $\La_1$ is the constant $n$-tuple ($\mu_1, \dots,
\mu_1$) and $\La_t = \La_M$ is the multitype weight. For $1 < j <
t$, the weight $\La_j = ( \la_1^j, \dots, \la_n^j)$ is defined by
$\la_i^j = \mu_i$ for $ i \leq k_{j-1}$, and $\la_i^j =
\mu_{k_{j-1}+1}$ for $i > k_{j-1}.$

If (\ref{1}) is  the defining equation in some
multitype coordinates, we define a model hypersurface associated
to $M$ at $p$ as
\begin{equation} M_H = \{(z,w) \in \mathbb C^{n+1}\ | \
 v  = P \zz \}. \label{2.10}\end{equation}
In order to analyse   biholomorphisms between models,
we will use the following terminology. 
\\[2mm]
\defi{Let $\La=(\la_1, \dots, \la_n)$ be a
weight.
  A transformation
$$ w^* = w + g(z_1, \dots z_n, w),\ \ \ \ \ z_i^* = z_i
+ f_i(z_1, \dots z_n, w)$$ preserving form (\ref{vp}) is called

-- $\La$-homogeneous if $f_i$ is a $\Lambda$-homogeneous
polynomial of weighted degree $\la_i$ and
 $ g $
is a $\Lambda$-homogeneous polynomial of weighted degree one,

-- $\La$-subhomogeneous if $f_i$ is a polynomial consisting of
monomials of weighted degree less or equal to $\la_i$ and $g$
consists of monomials of weighted degree less or equal to one,

-- $\La$-superhomogeneous if the Taylor expansion of $f_i$
consists of terms of weighted degree greater or equal to $\la_i$
and $g$
 consists of
terms of weighted degree  greater or equal to one.
 }
\\[2mm]

Note that we only consider nonsingular transformations (with
nonvanishing Jacobian at the origin).

 We now fix $\La_M$-adapted coordinates, and write the
corresponding leading polynomial $P$ as

\begin{equation}
P(z, \bar z) = \sum_{|(\alpha, \hat \alpha)|_{\La_M} = 1} A_{\al,
\hat \al} z^{\al} \bar z^{\hat \al}. \label{pal}
\end{equation}

 Let $P^k$ denote the restriction of $P$ to the first $k$
coordinate axes,

$$ P^k(z_1, \dots, z_k, \bar z_1, \dots, \bar z_k)
= P(z_1, \dots, z_k, 0, \dots, 0, \bar z_1, \dots, \bar z_k, 0, \dots, 0).$$

 It follows from the definition that  $\La_M$-homogeneous transformations are of the form

\begin{equation}
z^*_i = z_i +\sum_{|\al|_{\La_M} = \mu_i}C_{\al}z^{\al}, \ \ \ \
\ w^* = c w + \sum_{|\al|_{\La_M} = 1}D_{\al}z^{\al}
\end{equation}
where
$c \in
\mathbb R^*$.

The set of  such transformations forms a group, which will be
denoted by ${\mathcal H}$.
The subgroup of ${\mathcal H}$ consisting of transformations for
which $g=0$ (preserving the $w$ variable) will be denoted by
$\mathcal {H^Z}$.
Finally, let  $\mathcal L$ denote the subgroup of $\mathcal{H^Z}$,
consisting of all linear transformations in $\mathcal{H^Z}$.

\section{A normalization of the model}

We will use  the truncated leading polynomial $P^k$, $k
=1, \dots, n$,
to define a normalization condition corresponding to $\La_M$-homogeneous 
changes in the $z_k$ variable.
\\[3mm]
\defi{ Multitype  coordinates $(z,w)$   are called regular, if
for each $k = 1, \dots, n,$
\begin{equation}\frac{\partial P^k}{\partial z_k} (z_1, \dots, z_k, \bar z_1, \dots
, \bar z_k) \label{pk}
\end{equation}
 is not identically zero.}
\\[2mm]
The following
lemma shows that  regular coordinates do exist, and are in fact
generic among multitype coordinates.
\\[2mm]
 \lema{Let (z,w) be multitype coordinates. Then
 there exists a transformation $H \in \mathcal{H^Z}$,
 such that the new coordinates are regular.
 }
 \\[2mm]
 {\it Proof:} The proof is by induction. We will assume that
$\frac{\partial P^j}{\partial z_j} $
 is not identically zero
 for all $j < k$, and find  transformations
 which preserve this and attain the condition  for $P^k$.

 Let $k'$ be the largest index such that $\mu_k = \mu_{k'}$.
 Clearly, $P^{k'}$ has to depend on $z_k$, otherwise
 we could
 lower the weight of $z_k$ and obtain a lexicographically smaller
 distinguished weight, contradicting the definition of
 $\La_M$.
 Pick any monomial in
 $P^{k'}$ containing $z_k$, say
$$A_{\be, \hat \be}z^{\be}\bar z^{\hat \be},$$
where $A_{\be, \hat \be} \neq 0$, and $(\be, \hat \be)$ satisfies
$\be_j = \hat \be_j =0$ for $j>k'$, and $\be_k + \hat \be_k \neq
0$.
  Consider all terms in $P^{k'}$ with the same
initial part in the   variables $z_1, \dots, z_{k-1}$,
  $$ \left(\prod_{j<k}
   z_j^{\be_j}\bar z_j^{\hat \be_j}\right)Q(z_k, \dots, z_{k'}, \bar z_k, \dots, \bar z_{k'}).$$
  Clearly, for  a generic  linear transformation of the
  variables $z_k, \dots, z_{k'}$, in the new coordinates
  the corresponding homogeneous polynomial $Q^*$ does not
  vanish on the $z_k$ axis. It follows that
the restriction of  $P^*$ to $z_{k+1} = \dots = z_n = 0$ depends
on $z_k$. This finishes the proof.
\\[2mm]
 The following definition singles out a leading term in $P$ for
 each of the variables.
\\[2mm]

\defi{Let (z,w) be regular
coordinates. The leading term in the variable $z_k$ is given by the
lexicographically smallest multiindex pair $\Gamma^k = (\gamma^k,
\hat \gamma^k)$, such that
\begin{equation}
\gamma^k_j = \hat \gamma^k_j = 0 \ \ \ \text{for} \ \ \  j = k+1,
\dots, n,
\end{equation}

\begin{equation}
\ga_k^k + \hat \ga_k^k \neq 0 \ \ \ \text{and}\ \ \ \ \
A_{\gamma^k, \hat \gamma^k} \neq 0.
\end{equation}
}
\\[2mm]
The leading terms
are
 used to define a
normalization of $P$.
\\[2mm]
\defi{
Let $(z,w)$ be regular coordinates.
 We say that the leading polynomial $P$, given by (\ref{pal}),  is
normalized if for every $k$
$$ (i)\ \ \ \ \ \ \ \ A_{\ga^k, \hat \ga^k} = 1$$
and
$$(ii)\ \ \ \ \ \ \ \ \ A_{\al, \hat \al} = 0$$ for any multiindex pair $(\al, \hat \al)$ such
that $ \al = \ga^k$, $\ \hat \al_j = \hat \gamma^k_j  \ \text{for}
\ j < k,$, $ \hat \al_k = \ga_k^k -1$  and $\vert \hat \al\vert_{\La_M} = \vert \hat
\ga^k\vert_{\La_M}.$ }
\\[2mm]
 We will denote by $\ep^k$ the multiindex of length $n$ whose
$k$-th component is equal to one and other components are zero.

It is straightforward to show that the  normalization of $P$ can
indeed be attained by a $\La_M$-homogeneous transformation.
\\[2mm]
\lema{There exist regular coordinates in which the leading
polynomial $P(z, \bar z)$ is
normalized. }
\\[2mm]
{\it Proof:} By induction. Let us assume that we have found regular
coordinates such that (i) and (ii) are satisfied for all
$\Gamma^j$ with $j<k$ (note that $\Gamma^j$ are determined by the
coordinates). We will change the variable $z_k$ in such a way that
(i) and (ii) is satisfied also for $\Gamma^k$.
 The
transformation will be of the form
\begin{equation}
z_k  =
\sum_{\vert \al \vert_{\La_M} = \mu_k }C_{\al} (z^*)^{\al}, \ \ \
\ \ z_j = z_j^* \ \ {for}\ \  j\neq k,
\label{aab}
\end{equation} 
where $C_{\alpha}\neq 0$ implies $\al_j = 0$ for $j < k.$
 Substituting into $v = P(z, \bar z)$, we determine
the coefficients which attain the normalization condition for $k$.
This gives
$$ C_{\al} = - A_{\ga^k, \hat \ga^k - \ep^k + \al} + \dots,$$
for $\al \neq \ep^k$.
The condition $A_{\ga^k, \hat \ga^k} = 1$ is attained by taking
 $C_{\ep^k}$ as  a solution to
$$A_{\ga^k, \hat \ga^k}(C_{\ep^k})^{\ga_k^k}
 (\bar C_{\ep^k})^{\hat \ga_k^k} = 1.$$
That finishes the proof.
\\[2mm]

\section{Biholomorphic equivalence of models}
In this section we consider the local equivalence problem for models.
We start by showing that a transformation
preserves form (\ref{1}) if and only
if it is superhomogeneous.
\\[2mm]
 \prop{  \ A biholomorphic transformation
takes $\La_M$-adapted
coordinates into
$\La_M$-adapted
coordinates if and only if it is
$\La_M$-superhomogeneous. }
\\[2mm]
 {\it Proof.\ }
We first prove the only if part of the statement.
 Consider a transformation
\begin{equation}\begin{aligned}
  z^* = z + & f(z,w)\\
  w^* = w + & g(z,w), \\
\end{aligned} \label{fg} \end{equation}
which takes $\La_M$-adapted coordinates $(z,w)$ into $\La_M$-adapted coordinates $(z^*, w^*)$.

Let $v^* =F^*(z^*, \bar z^*, u^*)$ be the defining equation of $M$ in the
new coordinates.
 Substituting (\ref{fg}) into $v^* = F^*(z^*, \bar z^*, u^*)$
 we obtain the transformation formula
\begin{equation}\begin{aligned}  F^*(z + f(z,u+iF(z, \bar z, u)),
\overline{z + f(z,u+iF(z, \bar z, u))}, & u \;+ \\ +\;Re\
 g
 (z,u+iF(z, \bar z, u))   = F(z, \bar z, u) +
  Im\ g(z, u+&iF(z, \bar z,
u)).\label{covf}\end{aligned} \end{equation}

 Without any loss of generality,
we may assume that $P$ is
normalized (applying a $\La_M$-homogeneous transformation in the
source space, if necessary). On the other hand, 
 we do not assume that $P^*$
is normalized. Instead, in the target space we use an element of
$\mathcal{L}$ to normalize the linear part of the transformation
and assume that the Jacobi matrix of the transformation at the
origin is the unit matrix.

 By induction we will show that the transformation has to be
superhomogeneous with respect to all  weights in the generating
sequence
 $\La_1,
\La_2, \dots, \La_t$.

For $l=1$, we have $\La_1 = (m_1, ..., m_1)$. Hence $\La_1$-homogeneous
 transformations in $\mathcal{H^Z}$ are linear,
and the claim is obvious.

Let $l>1$ and
assume the transformation is $\La_j$-superhomogeneous for all
$j$ with $j<l$. We will prove that it is also
$\La_l$-superhomogeneous.
 Note that $\la^l_j < \la^{l-1}_j$
 if and only if $j > k_{l-1}.$

We separate the strictly subhomogeneous part (with respect to $\La_l$) of the inverse transformation,
and write
\begin{equation}
z_i = z^*_i + \sum_{\vert \alpha\vert_{\La_{l}} < \la^l_i}
C^i_{
\alpha}
(z^*)^{\alpha} + O_{\La_l}(\la_i^l). \label{11}
\end{equation}
Note that $w = w^* + o_{\La_l}(1)$, since $P$ and $P^*$ contain no
pluriharmonic terms.
In this notation, let
 $$ \Theta = \{ (i,\alpha) \in \mathbb Z^{n+1};\  C^i_{\al} \neq 0\}.$$
The elements of $\Theta$ have the following immediate properties.
 If $(i,\alpha)
\in \Theta$, then $\al_j = 0 $ for $j\leq i$, since $\la_j^l \geq
\la_i^l$.
Further, by $\La_{l-1}$-superhomogeneity, each of the terms
appearing in (\ref{11}) must contain at least one of the variables
$z_{k_{l-1} +1}, \dots, z_n$. We denote
$$ S(\al) = \al_{k_{l-1}+1} + \dots + \al_n.$$
Hence $S(\al) > 0$ for all $(i,\alpha) \in \Theta$.

Analogous notation will be also used for multiindex pairs:
$$S(\al, \hat \al ) =\al_{k_{l-1}+1} + \dots + \al_n + \hat \al_{k_{l-1}+1} +
\dots +\hat \al_n.
$$

Let $m_S$ be the minimal value of $S(\al)$ as $(i,\al)$ ranges
over $\Theta$.
For $(i,\alpha) \in \Theta$ consider the "gap"
$$G(i,{\al}) = \la^l_i -\sum_{j=1}^n \al_j
\la^l_j.$$

Among all pairs $(i,\al)$ in $\Theta $ for which $S(\al) = m_S$, let $\Xi$
denote the set of those for which $G(i,\al)$ is maximal.
 Next, let $m$ be the smallest integer
such that $(m,\al) \in \Xi$ for some $\al$. Now we fix one such
pair, $(m,\delta)\in \Xi$ and consider the corresponding monomial in \ref{11}:

$$ C^m_{(0,\dots,0,\de_{m+1},
\dots, \de_n) } \prod_{j > m} (z^*_j)^{\de_j},$$ where
  $\de = (0,\dots,0,\de_{m+1}
\dots, \de_n).$ Note that $ m \leq k_{l-1}$.

Substituting  (\ref{11}) into
$$v =\sum_{|(\al, \hat \al)|_{\La_l}=1} A_{\al, \hat \al } z^{\al} \bar z^{\hat \al} + o_{\La_l}(1), $$
we  compute the coefficient of
\begin{equation}
(z^*)^{\ga^m }(\bar z^*)^{\hat \ga^m - \ep^m + \de}. \label{zgd}
\end{equation}
Since $F$ starts with weight one, it is enough to consider the
strictly subhomogeneous part of the transformation. Hence we need
to consider the expansion of
$$F(z^*_1 + \sum_{\vert \alpha\vert_{\La_{l}} < \la^l_1}
 C^1_{\al}(z^*)^{\al}, \dots, z^*_n + \sum_{\vert \alpha
 \vert_{\La_{l}} < \la^l_
n}
  C^n_{\al}(z^*)^{\al},
\overline{z_1^* + 
\dots, \phantom{}}, 0)$$ 
First, consider terms
coming from
the leading polynomial.
If for some
multiindex pair $(\be, \hat \be)$ the coefficient
 $A_{\be, \hat
\be}$ enters the equation for (\ref{zgd}), 
then by the choice of $(m, \delta)$ 
there exists a multiindex $\al$ and $ j\in \{ 1, \dots, n\}$
such that
$$\hat \ga^m - \ep^m + \de = \hat \be - \ep^j + \al,$$
and
$\be = \ga^m$. But, again  by the choice of $(m, \de)$, the
gaps satisfy
$$ |\de - \ep^m|_{\La_l}  = |\al - \ep^j|_{\La_l},$$
so
$$|\hat \ga^m|_{\La_l} = |\hat \be|_{\La_l}.$$
Moreover, $\hat \ga^m_j = \hat \be_j$ for all $j<m$, which gives
contradiction with the normalization of $P$.
Note that 
$$ \ga^m -\ep_m
=  \be -\ep_j$$ 
is impossible, since it forces $\la_j = \la_m$,
 which contradicts the linear part of the
transformation being the identity.

It remains to prove that  terms of weight greater than one in $F$
do not enter the equation for (\ref{zgd}).
Let $F_{\al, \hat \al, l} z^{\al} \bar z^{\hat \al } u^l $ be such a term, where
$\vert (\al, \hat \al) \vert_{\La_M} > 1$. 
By the choice  of $(m,\de)$, in order to
influence a term of weight $1-G(m,\de)$ in $F^*$
we have  to substitute
 at least twice a term with the lowest
value of $S(\al)$,
or a term with a higher value of $S(\al)$.
In both cases the  value of $S(\al, \hat \al)$
for the resulting term  is bigger  than
 $$m_S = S(\ga^m - \ep^m + \de, \hat \ga^m).$$ Thus we
 have proved that a transformation
 which takes
$\La_M$-adapted
coordinates into $\La_M$-adapted
coordinates  is $\La_M$-superhomogeneous.
 The converse follows immediately from (\ref{covf}).
\\[2mm]

 Now we can describe explicitly biholomorphisms between different
models.
\\[2mm]
\prop{Let $M_H$ and $\tilde M_H$ be two models for $M$ at $p$.
Then there is a $\La_M$-homogeneous transformation which maps $M_H$ to
$\tilde M_H$. In particular, all models are biholomorphically
equivalent by a polynomial transformation.}
{\it Proof: } By the previous proposition, the coordinates in which
$M_H$ is the model are related to those in which $\tilde M_H$ is
the model by a  $\La_M$-superhomogeneous transformation.
But terms of weight greater than $\la_i$ in $f_i$ influence only
terms of weight greater then one in $F^*$. Hence $\tilde M_H$ is
obtained by the homogeneous part of this transformation.

\section{Computing the multitype}

Using Theorem 4.1,  the process of computing multitype
can be described as follows.

In the first step, we consider local holomorphic coordinates in which the leading
polynomial in the variables $z, \bar z $ contains no pluriharmonic term. The first multitype
component $m_1$ is then equal to
the degree of this polynomial. Hence $m_1 = \frac1{\mu_1}$ is equal to the Bloom-Graham type
of $M$ at $p$,
 and we set $\La_1$ =
($\mu_1, \dots, \mu_1$).

 In the second step, consider all
$\La_1$-homogeneous transformations 
and choose  one
which makes the leading polynomial $P_1$ independent of the largest
number of variables. Let $d_1$ denote this number. Permuting
variables, if necessary, we can assume that in such coordinates,
 $$v = P_1(z_1, \dots, z_{n-d_1} , \bar z_1, \dots \bar z_{n-d_1} ) + Q_1(z, \bar
 z) + o(u),$$
 where $P_1$ is $\La_1$-homogeneous of weighted degree one, and $Q_1$ is $o_{\La_1}(1)$.
Since $\La_1$-homogeneous transformations are linear, 
and using the fact that
 for any weight $\La$ which is lexicographically smaller
 than $\La_1$, $\La$-adapted coordinates are also $\La_1$-adapted, it follows that
 $\mu_1 = \mu_2 = \dots = \mu_{n-d_1}$ and $\mu_{n-d_1 +1} < \mu_1$.
Let $$Q_1(z, \bar z) =\sum_{\vert (\al, \hat \al) \vert_{\La_1} > 1}
C^1_{\al, \hat \al} z^{\al} \bar z^{\hat \al},$$ and denote
$$\Theta_1 = \{ (\al, \hat \al) \ |\ C^1_{\al, \hat \al} \neq 0 \ \text{and}\
\sum_{i=1}^{n-d_1} \al_i + \hat \al_i < m_1 \}.
$$
For each 
$(\be, \hat \be) \in \Theta_1$
 consider the number
\begin{equation} W_1(\be, \hat \be) =
\frac{1 - \sum_{i=1}^{n-d_1} (\be_i + \hat \be_i)\mu_1}{
  \sum_{i = n-d_1+1}^n \be_i + \hat \be_i}.
\label{ww}
  \end{equation}
The weight $\La_2$ is defined by  $\la^2_j = \mu_1$ for $j \leq
n-d_1$, and
  $$ \la_j^2 = \max_{(\al, \hat \al) \in \Theta_1} W_1(\al, \hat \al)$$
for $j > n-d_1$.

By the definition of $\La_2$, the leading polynomial with respect
to $\La_2$ in the above coordinates depends on at least $n-d_1 +1$
variables. This ends the second step. 

 Now we continue the process.
 In the j-th step, $j>2$,  we use the coordinates obtained in the previous step, 
and consider all $\La_{j-1} $-homogeneous transformations.
We denote by $d_{j-1}$ the largest number of variables which do not appear in the
 leading polynomial after the transformation, and 
 fix such a coordinate system.
It is easy to show, using the same arguments as in Theorem 4.1., that
any transformation which takes $\La_{j-1}$-adapted coordinates into  $\La_{j-1}$-adapted
coordinates has to be $\La_{j-1}$-superhomogeneous.
 If $d_{j-1} < d_{j-2}$, using this 
and  the fact that
 for any weight $\La$ which is lexicographically smaller
 than $\La_{j-1}$, $\La$-adapted coordinates are also $\La_{j-1}$-adapted,
it follows that we have
 determined the $(d_{j-2} - d_{j-1})$ multitype entries
 $$ \mu_{n-d_{j-2}+1} = \dots = \mu_{n-d_{j-1}} = \la^{j-1}_{n-d_{j-2}+1},$$
and set $\la_i^j = \mu_i$ for $i\leq n-d_{j-1}.$
To define the remaining entries of $\La_j$, we write
 $$v = P_{j-1}(z_1, \dots, z_{n-d_{j-1}} , \bar z_1, \dots \bar z_{n-d_{j-1}} ) + Q_{j-1}(z, \bar
 z) + o(u),$$
 where $P_{j-1} $ is $\La_{j-1}$-homogeneous of weighted degree one, and $Q_{j-1}$ is $o_{\La_{j-1}}(1)$, 
$$Q_{j-1}(z, \bar z) =\sum_{\vert (\al, \hat \al)
\vert_{\La_{j-1}}
> 1} C^{j-1}_{\al, \hat \al} z^{\al} \bar z^{\hat \al}.$$
Let
$$\Theta_{j-1} = \{ (\al, \hat \al) \ |\ C^{j-1}_{\al, \hat \al} \neq 0 \ \text{and}\
\sum_{i=1}^{n-d_{j-1}} (\al_i + \hat \al_i)\mu_i < 1 \}.
$$
As before, denote
\begin{equation} W_{j-1}(\be, \hat \be) =
\frac{1 - \sum_{i=1}^{n-d_{j-1}} (\be_i + \hat \be_i)\mu_i}{
  \sum_{i = n-d_{j-1}+1}^n \be_i + \hat \be_i}.
\label{ww2}
  \end{equation}
The remaining entries of  $\La_j$ are defined by
  $$ \la_i^j = \max_{\al \in \Theta} W_{j-1}(\al, \hat \al)$$
for $j > n-d_{j-1}$.

If $d_{j-1} = d_{j-2}$, we only use (\ref{ww2}) to determine
$\la^j_{n-d_{j-1} +1}, \dots \la^j_n$. No multitype component is
determined at this step. 

It is immediate to verify that the
process terminates after finitely many steps,
 and  determines all components of the multitype
 weight.

\end{document}